\documentclass[12pt]{amsart}
\usepackage{amscd}
\usepackage{amsmath}
\usepackage{amsxtra}
\usepackage{amsfonts}
\usepackage{amssymb}

\newtheorem{theorem}{Theorem}[section]
\newtheorem{corollary}[theorem]{Corollary}
\newtheorem{lemma}[theorem]{Lemma}
\newtheorem{proposition}[theorem]{Proposition}

\theoremstyle{definition}
\newtheorem{definition}[theorem]{Definition}
\newtheorem{remark}[theorem]{Remark}

\theoremstyle{remark}

\renewcommand{\theclaim}{\textup{\theclaim}}

\numberwithin{equation}{section}

\def\openone

{\mathchoice

{\hbox{\upshape \small1\kern-3.3pt\normalsize1}}

{\hbox{\upshape \small1\kern-3.3pt\normalsize1}}

{\hbox{\upshape \tiny1\kern-2.3pt\SMALL1}}

{\hbox{\upshape \Tiny1\kern-2pt\tiny1}}}

\makeatletter

\newbox\ipbox

\newcommand{\ip}[2]{\left\langle #1\, , \,#2\right\rangle}
\newcommand{\diracb}[1]{\left\langle #1\mathrel{\mathchoice

{\setbox\ipbox=\hbox{$\displaystyle \left\langle\mathstrut
#1\right.$}

\vrule height\ht\ipbox width0.25pt depth\dp\ipbox}

{\setbox\ipbox=\hbox{$\textstyle \left\langle\mathstrut
#1\right.$}

\vrule height\ht\ipbox width0.25pt depth\dp\ipbox}

{\setbox\ipbox=\hbox{$\scriptstyle \left\langle\mathstrut
#1\right.$}

\vrule height\ht\ipbox width0.25pt depth\dp\ipbox}

{\setbox\ipbox=\hbox{$\scriptscriptstyle \left\langle\mathstrut
#1\right.$}

\vrule height\ht\ipbox width0.25pt depth\dp\ipbox}

}\right. }

\newcommand{\dirack}[1]{\left. \mathrel{\mathchoice

{\setbox\ipbox=\hbox{$\displaystyle \left.\mathstrut
#1\right\rangle$}

\vrule height\ht\ipbox width0.25pt depth\dp\ipbox}

{\setbox\ipbox=\hbox{$\textstyle \left.\mathstrut
#1\right\rangle$}

\vrule height\ht\ipbox width0.25pt depth\dp\ipbox}

{\setbox\ipbox=\hbox{$\scriptstyle \left.\mathstrut
#1\right\rangle$}

\vrule height\ht\ipbox width0.25pt depth\dp\ipbox}

{\setbox\ipbox=\hbox{$\scriptscriptstyle \left.\mathstrut
#1\right\rangle$}

\vrule height\ht\ipbox width0.25pt depth\dp\ipbox}

} #1\right\rangle}

\newcommand{\cj}[1]{\overline{#1}}

\newcommand{\bz}{\mathbb{Z}}

\newcommand{\bc}{\mathbb{C}}
\newcommand{\bt}{\mathbb{T}}
\newcommand{\bn}{\mathbb{N}}
\newcommand{\tr}{\textup{tr}}

\def\blfootnote{\xdef\@thefnmark{}\@footnotetext}

\input xy
\xyoption{all}
\usepackage{amssymb}

\def\N{\mathbb{N}}

\def\H{\mathcal{H}}

\def\-{^{-1}}

\def\D{\mathcal{D}}

\def\U{\mathcal{U}}

\def\K{\mathcal{K}}

\begin{document}

\title[Parseval frames for ICC groups]{Parseval frames for ICC groups}
\author{Dorin Ervin Dutkay}
\blfootnote{Research supported in part by a grant from the
National Science Foundation DMS-0704191}
\address{[Dorin Ervin Dutkay] University of Central Florida\\
    Department of Mathematics\\
    4000 Central Florida Blvd.\\
    P.O. Box 161364\\
    Orlando, FL 32816-1364\\
U.S.A.\\} \email{ddutkay@mail.ucf.edu}
\author{Deguang Han}
\address{[Deguang Han] University of Central Florida\\
    Department of Mathematics\\
    4000 Central Florida Blvd.\\
    P.O. Box 161364\\
    Orlando, FL 32816-1364\\
U.S.A.\\} \email{dhan@pegasus.cc.ucf.edu}
\author{Gabriel Picioroaga}
\address{[Gabriel Picioroaga]Binghamton University\\
         Department of Mathematical Sciences\\
         Binghamton, NY, 13902-6000\\
U.S.A.\\} \email{gabriel@math.binghamton.edu}
\thanks{}

\begin{abstract}
   We analyze Parseval frames generated by the action of an
   ICC group on a Hilbert space. We parametrize the set of all
   such Parseval frames by operators in the commutant of the
   corresponding representation. We characterize when two such frames
   are strongly disjoint. We prove an undersampling result showing that
   if the representation has a Parseval frame of equal norm vectors of norm $\frac{1}{\sqrt{N}}$,
   the Hilbert space is spanned by an orthonormal basis generated by a
   subgroup. As applications we obtain some sufficient conditions
   under which a unitary representation admits a Parseval frame
   which is spanned by an Riesz sequences generated by a
   subgroup. In particular, every subrepresentation of the left
   regular representation of a free group has this property.
\end{abstract}
\maketitle \tableofcontents

\section{Introduction}

Frames play a fundamental role in signal processing, image and
data compression and sampling theory. They provide an alternative
to orthonormal bases, and have the advantage of possessing a
certain degree of redundancy which can be useful in applications, for example
when data is lost during transmission. Also, frames can be better
localized, a feature which lead to the success of Gabor frames and
wavelet theory (see e.g., \cite{Dau92}).

 The term ``frame'' was introduced by Duffin and Schaffer
 \cite{DuSc} in their study of non-harmonic Fourier series, and
 has generated important research areas and remarkable
 breakthroughs \cite{Sei}. Recent results show that frames can
 provide a universal language in which many fundamental problems
 in pure mathematics can be formulated: the Kadison-Singer problem
 in operator algebras, the Bourgain-Tzafriri conjecture in Banach
 space theory, paving T\" oplitz operators in harmonic analysis
 and many others (see \cite{Cas06} for an excellent account).

\begin{definition}\label{def1}
Let $\H$ be a Hilbert space. A family of vectors $\{e_i\,|\,i\in
I\}$ in $\H$ is called a {\it frame} if there exist constants
$A,B>0$ such that for all $f\in \H$,
$$A\|f\|^2\leq\sum_{i\in I}|\ip{f}{e_i}|^2\leq B\|f\|^2.$$
If $A=B=1$, then $\{e_i\,|\, i\in I\}$ is called a {\it Parseval
frame}. If we only require the right hand inequality, then we say
that $\{e_i\,|\,i\in I\}$  is a {\it Bessel sequence}.
\end{definition}

In their seminal paper by Han and Larson \cite{HaLa00}, operator
theoretic foundations  for frame theory and group representations
where formulated. One of the key observations in their paper was
that every Parseval frame is the orthogonal projection of an
orthonormal basis. This lead them to the notion of {\it
disjointness} of frames.

\begin{definition}\label{def2}
Let $\{e_i\,|\, i\in I\}$, $\{f_i\,|\, i\in I\}$ be two Parseval
frames for the Hilbert spaces $\H_1$ and $\H_2$ respectively. The
Parseval frames are called {\it strongly disjoint} if $\{e_i\oplus
f_i\,|\, i\in I\}$ is a Parseval frame for $\H_1\oplus \H_2$.
Similarly, if we have $\{e_i^1\,|\,i\in I\},\dots,\{e_i^N\,|\,i\in
I\}$ Parseval frames for the Hilbert spaces $\H_i$ respectively,
we say that these frames form an {\it $N$-tuple of strongly
disjoint Parseval frames} if $\{e_i^1\oplus\dots\oplus
e_i^N\,|\,i\in I\}$ is a Parseval frame for
$\H_1\oplus\dots\oplus\H_N$.

 The Parseval frames are called {\it unitarily equivalent} if there exists a
 unitary operator $U:\H_1\rightarrow\H_2$ such that $Ue_i=f_i$ for all $i\in I$.
\end{definition}

If the direct sum of two Parseval frames is an orthonormal basis,
we say that one is the {\it complement} of the other, or one {\it
complements} the other.

Many properties of frames are encoded in the associated {\it frame
transform} (or {\it analysis operator}). It is the operator that
associates to a vector its coefficients in the given frame.

\begin{definition}\label{def3}
Let $\mathcal E:=\{e_i\,|\, i\in I\}$ be a Bessel sequence in a
Hilbert space $\H$. The operator $\Theta_{\mathcal
E}:\H\rightarrow l^2(I)$ defined by
$$\Theta_E(f)=(\ip{f}{e_i})_{i\in I}, \quad (f\in \H)$$
is called the {\it analysis operator or frame transform}
associated to $\mathcal E$.
\end{definition}

Strong disjointness can be characterized in terms of frame
transforms as follows:

\begin{proposition}\label{pr4}\cite{HaLa00}
Let $\mathcal E:=\{e_i\,|\, i\in I\}$ and $\mathcal
F:=\{f_i\,|\,i\in I\}$ be two Parseval frames for the Hilbert
spaces $\H_1$ and $\H_2$ respectively. Then $\mathcal E$ and
$\mathcal F$ are unitarily equivalent iff the frame transforms
$\Theta_{\mathcal E}$ and $\Theta_{\mathcal F}$ have the same
range. The following affirmations are equivalent:
\begin{enumerate}
    \item The Parseval frames $\mathcal E$ and $\mathcal F$ are strongly disjoint.
    \item The frame transforms $\Theta_{\mathcal E}$ and $\Theta_{\mathcal F}$ have orthogonal ranges.
    \item $\Theta_{\mathcal F}^*\Theta_{\mathcal E}=0$.
    \item For all $v_1\in \D_1$, $v_2\in \D_2$, where $\D_1$, $D_2$ are dense in $\H_1,\H_2$ respectively,
    $$\sum_{i\in I}\ip{v_1}{e_i}\cj{\ip{v_2}{f_i}}=0.$$
\end{enumerate}
\end{proposition}
In what follows we will call two Bessel sequences (with same index
set) {\it strongly disjoint} if the range spaces of their analysis
operators are orthogonal.

In this paper we will be interested only in the Parseval frames
generated by the action of an infinite-conjugacy-classes (ICC)
group on a Hilbert space. The ICC property implies that the
associated left-regular representation is a $II_1$ factor, and we
can use the theory of factor von Neumann algebras \cite{JoSu97}.
Our purpose is to follow and further the results in \cite{HaLa00}.

After the Han-Larson paper, frames for abstract Abelian groups
have been studied in \cite{TaWe,Web} using Pontryagin duality. As
we shall see, Parseval frames for ICC groups have quite different
properties. Many strongly disjoint Parseval frames can be found in
the same representation (Theorem \ref{thi3}), and there are
undersampling results in some cases (Proposition \ref{prp1} and
Theorem \ref{thp2}). Frames for ICC groups fit into the more
general theory of group-like unitary systems, a theory which has
the Gabor (or Weyl-Heisenberg) frames as the central example, see
\cite{GaHa03,HL01,Han3} for details.

\begin{definition}\label{def6}
Let $G$ be a countable group. Let $\pi$ be a unitary
representation of $G$ on the Hilbert space $\H$. A vector $\xi\in
\H$ is called a {\it frame/ Parseval frame/ ONB vector} for $\H$
(with the representation $\pi$) iff $\{\pi(g)\xi\,|\,g\in G\}$ is
a frame/Parseval frame/ONB for $\H$.

Two Parseval frame vectors $\xi, \eta$ for $\H$ are called {\it
unitarily equivalent/strongly disjoint} if the corresponding
Parseval frames $\{\pi(g)\xi\,|\,g\in G\}$ and $\{\pi(g)\eta\,|\,
g\in G\}$ are unitarily equivalent/ strongly disjoint. Similarly
for an {\it $N$-tuple of strongly disjoint Parseval frame
vectors}.
\end{definition}

As shown in \cite{HaLa00}, every representation of a group that
has a Parseval frame vector is isomorphic to a sub-representation
of the left-regular representation (Proposition  \ref{propi1}).

\begin{definition}\label{def5}
Let $G$ be a countable group. The {\it left-regular
representation} $\lambda$ of $G$ is defined on $l^2(G)$ by
$$(\lambda(g)\xi)(h)=\xi(g^{-1}h),\quad(\xi\in l^2(G), h\in G, g\in G).$$
Equivalently, if $\delta_g$, $g\in G$ is the canonical orthonormal
basis for $l^2(G)$, then
$$\lambda(g)\delta_h=\delta_{gh},\quad(h,g\in G).$$
The von Neumann algebra generated by the unitary operators
$\lambda(g)$, $g\in G$ is denoted by $L(G)$.

The {\it right-regular representation} $\rho$ of $G$ is defined on
$l^2(G)$ by
$$(\rho(g)\xi)(h)=\xi(hg),\quad(\xi\in l^2(g), h\in G, g\in G).$$
Equivalently,
$$\rho(g)\delta_h=\delta_{hg^{-1}},\quad(h,g\in G).$$

The group $G$ is called an ICC (infinite conjugacy classes) group
if, for all $g\in G$, $g\neq e$, the set $\{hgh^{-1}\,|\, h\in
G\}$ is infinite.
\end{definition}

The commutant of the von Neumann algebra $L(G)$ is the von Neumann
algebra generated by the right-regular representation $\rho$. When
the group is ICC, $L(G)$ is a $II_1$ factor (see \cite{JoSu97}).

The paper is organized as follows: in Section \ref{gene}, we study
the general properties of Parseval frames for ICC groups. In
Theorem \ref{thi3} we show how a Parseval frame vector can be
complemented by several other Parseval frame vectors for the same
representation and a ``remainder'' Parseval frame vector for a
subspace. While Lemma \ref{pri4} characterizes strongly disjoint
Parseval frame vectors in terms of their cyclic projections we
present some new properties along the same lines in Proposition
\ref{pri5}.

We give a parametrization of Parseval frame vectors in Theorem
\ref{pri6}. Another such parametrization was given in
\cite{HaLa00}, see Theorem \ref{thi2}. The advantage of our
parametrization is that it uses operators in the commutant in the
representation, so it can be extended to {\it all} the vectors of
the frame. Theorem \ref{pri7} characterizes strong disjointness
and unitary equivalence in terms of this parametrization.

In Section \ref{pars} we will be interested in the relation
between Parseval frames and subgroups. We will work with Parseval
frame vectors that have norm-square equal to $\frac{1}{N}$ with
$N\in\bn$ and we will assume that there exists a subgroup $H$ of
index $N$. We give two ``undersampling'' results: in Proposition
\ref{prp1}, we show that in this situation, there exist
orthonormal bases generated by the action of the subgroup $H$, and
in Theorem \ref{thp2} we show that, if in addition $H$ is normal,
then we can construct $N$ strongly disjoint Parseval frame vectors
such that, when undersampled to the subgroup, they form
orthonormal bases (up to a multiplicative constant). We apply
these results to the Feichtinger frame decomposition problem for
free groups to get that every frame representation for free groups
admits a frame that is a finite union of Riesz sequences.

\begin{definition}\label{def7}
We will use the following notations: if $\mathcal A$ is a set of
operators on a Hilbert space $\H$, then $\mathcal A'$ is the
commutant of $\mathcal A$, i.e., the set of all operators that
commute with all the operators in $\mathcal A$. $\mathcal A''$ is
the double-commutant, i.e., the commutant of $\mathcal A'$. By von
Neumann's double commutant theorem, $\mathcal A''$ coincides with
the von Neumann algebra generated by $\mathcal A$. Two
(orthogonal) projections $p$ and $q$ in a von Neumann algebra
$\mathcal A$ are said to be equivalent (denoted by $p\sim q$) if
there exists an operator (partial isometry) $u\in \mathcal A$ such
that $uu^{*} = p$ and $u^{*}u = q$. A von Neumann algebra
$\mathcal A$ is finite it there is no proper projection of $I$
that is equivalent to $I$ in $\mathcal A$. We refer to
\cite{KaRi97} for more details and some properties about von
Neumann algebras that will be used in the rest of the paper.

If $\pi$ is a representation of a group $G$ on a Hilbert space
$\H$, and $N\in\bn$, then $\pi^N$ is the representation of $G$ on
$\H^N$ defined by
$$\pi^N(g)=\underbrace{\pi(g)\oplus\dots\oplus\pi(g)}_{N\mbox{ times}},\quad(g\in G).$$

If $p'$ is a projection in $\pi(G)'$, then $\pi p'$ is the
representation of $G$ on $p'\H$ defined by
$$(\pi p')(g)=p'\pi(g)p',\quad(g\in G).$$

If $M$ is a von Neumann algebra with a trace, we will denote the
trace by $\tr_M$
\end{definition}

The next proposition proved by Han and Larson in \cite{HaLa00} is
the starting point for the theory of Parseval frames for groups.
It shows that any Parseval frame generated by the representation
of a group is in fact isomorphic to the projection of the
canonical basis in the left-regular representation of the group.
The isomorphism is in fact the frame transform, the projection is
its range and it lies in the commutant of the left-regular
representation.

\begin{proposition}\cite{HaLa00}\label{propi1}
Let $G$ be a countable ICC group and let $\pi:G\rightarrow\U(\H)$
be a unitary representation of the group $G$ on the Hilbert space
$\H$. Suppose $\xi\in\H$ is a Parseval frame vector for $\H$.
Then:
\begin{enumerate}
\item The frame transform $\Theta_\xi$ is an isometric isomorphism
between $\H$ and the subspace $p'l^2(G)$, where
$p':=\Theta_\xi\Theta_\xi^*$. \item The frame transform
$\Theta_\xi$ intertwines the representations $\pi$ on $\H$ and
$\lambda$ on $l^2(G)$, i.e.,
$\Theta_\xi\pi(g)=\lambda(g)\Theta_\xi$, for all $g\in G$. The
projection $p':=\Theta_\xi\Theta_\xi^*$ commutes with
$\lambda(G)$. \item $\Theta_\xi\xi=p'\delta_e$,
$\Theta_\xi^*\delta_e=\xi$. The trace of $p'$ in $L(G)'$ is
$\tr_{L(G)'}(p')=\|\xi\|^2$.
\end{enumerate}
Thus the Parseval frame $\{\pi(g)\xi\,|\,g\in G\}$ for $\H$ is
unitarily equivalent to the Parseval frame
$\{\lambda(g)p'\delta_e\,|\,g\in G\}$ for $p'l^2(G)$, via the
frame transform $\Theta_\xi$.
\end{proposition}

In \cite{HaLa00}, the authors proved that Parseval vectors can be
parametrized by unitary operators in the algebra $\pi(G)''$. Note
that these operators are {\it not} in the commutant, and this has
the drawback that one can not map the entire Parseval frame into
the other. We will give an alternative parametrization, that uses
operators in the commutant in Proposition \ref{pri6}.

\begin{theorem}\cite{HaLa00}\label{thi2}
Let $G$ be a countable ICC group and let $\pi:G\rightarrow\U(\H)$
be a unitary representation of the group $G$ on the Hilbert space
$\H$. Suppose $\xi\in\H$ is a Parseval frame vector for $\H$. Then
$\eta\in\H$ is a Parseval frame vector for $\H$ if and only if
there exists a unitary $u\in\pi(G)''$ such that $u\xi=\eta$. In
particular $\|\eta\|=\|\xi\|$.
\end{theorem}

\section{General theory: A new parametrization theorem}\label{gene}

We will consider an ICC group $G$ and $\pi:G\rightarrow\mathcal
U(\H)$ a unitary representation of the group on the Hilbert space
$\H$. We will assume that this representation has a Parseval frame
vector $\xi_1$. The main goal of this section is to obtain an
alternate parametrization of all the Parseval frame vectors by
using operator vectors with entries in the commutant of $\pi(G)$.
Although the focus of this paper is on ICC groups, the new
parametrization result works for other groups for which we will
discuss the details at the end of this section.

The norm-square of the vector $\|\xi_1\|^2$ fits into one $N$
times, and there may be some remainder $0\leq r\leq \|\xi_1\|^2$.
The next theorem shows that we can complement the Parseval frame
vector $\xi_1$ by $N-1$ Parseval frame vectors for $\H$ and one
``remainder'' Parseval frame vector for a subspace of $\H$.
Moreover, the complementing procedure works also if we start with
several strongly disjoint Parseval frame vectors for $\H$. Of
course, if $\|\xi_1\|^2=\frac{1}{N}$ with $N\in\bn$, the remainder
Parseval vector is not needed and can be discarded. Our results
have a simpler statement if the extra assumption
$\|\xi_1\|^2=\frac{1}{N}$ is added, and the remainders disappear.
We recommend the reader do this for an easier understanding of the
statements. However we sacrificed (part of) the aesthetics for
generality.

\begin{lemma} \cite{HaLa00, GaHa04} \label{lemma2.1}
Let $G$ be a countable ICC group and let $\pi:G\rightarrow\U(\H)$
be a unitary representation of the group $G$ on the Hilbert space
$\H$. Suppose $\xi\in\H$ is a Parseval frame vector for $\H$ and
let $p'_{\xi} = \Theta_{\xi}\Theta_{\xi}^{*}$ and  $p'$ be a
projection in $L(G)'$. Then

(i) $p'\sim p'_{\xi}$ in the von Neumann algebra $L(G)'$ if and
only if there exists a Parseval frame vector $\eta$ for $\H$ such
that $p' = p'_{\eta}$;

(ii) $p'$ is equivalent to a subprojecion of $p'_{\xi}$ if and
only if there exists a Parseval frame vector $\eta$ for
$\overline{span}\{\pi(g)\eta: g\in G\}$ such that $p' =
p'_{\eta}$;
\end{lemma}

\begin{theorem}\label{thi3}
Let $G$ be a countable ICC group and let $\pi:G\rightarrow\U(\H)$
be a unitary representation of the group $G$ on the Hilbert space
$\H$. Suppose $\xi_1\in\H$ is a Parseval frame vector for $\H$.
Let $N:=\left\lfloor \frac{1}{\|\xi_1\|^2}\right\rfloor \in\bn$.
\begin{enumerate}
\item There exist $\xi_2,\dots,\xi_N,\xi_{N+1}\in\H$ with the
following properties:
\begin{enumerate}
\item $\xi_2,\dots,\xi_N$ are Parseval frame vectors for $\H$.
\item $\|\xi_{N+1}\|^2=1-N\|\xi_1\|^2=:r<\|\xi_1\|^2$, and there
exists a projection $p'_{r}\in\pi(G)'$ such that
$\{\pi(g)\xi_{N+1}\,|\,g\in G\}$ is a Parseval frame for $p_r'\H$,
so $\xi_{N+1}$ is a Parseval frame vector for $p_r'\H$ with the
representation $\pi^r:=\pi p_r'$. \item $\xi_1,\dots,\xi_{N+1}$ is
a strongly disjoint $N+1$-tuple of Parseval frame vectors.
\end{enumerate}

\item If $\xi_1,\dots,\xi_{N+1}$ are as in (i) then there is no
vector $\xi_{N+2}$ such that $\xi_{N+2}$ is a Parseval frame
vector for some representation $\pi_{N+2}$ of $G$, and $\xi_{N+2}$
is strongly disjoint from all $\xi_i$, $i=1,N+1$.


\item If $1\leq M\leq N$ and $\xi_1,\dots,\xi_M$ is a strongly
disjoint $M$-tuple of Parseval frame vectors for $\H$, then there
exist $\xi_{M+1},\dots,\xi_N,\xi_{N+1}$ such that the properties
in (i) are satisfied.
\end{enumerate}
\end{theorem}

\begin{proof} (i) Let $p'_1 = \Theta_{\xi_{1}}\Theta_{\xi_{1}}^{*}$.  Since $L(G)'$ is a factor von Neumann algebra,
we have that there exist mutually orthogonal projections $p'_{2} ,
..., p'_{N}, p'_{N+1}$ in $L(G)'$ such that $p'_{i}\sim p'_{1}$
for $i =2, ..., N$,  $p'_{N+1}$ is equivalent to a subprojection
of $p'_1$ and $\sum_{i=1}^{N+1}p'_{i} = I$.  By Lemma
\ref{lemma2.1}, there exist $\xi_2\dots,\xi_N,\xi_{N+1}\in\H$ such
that
 $\xi_2,\dots,\xi_N$ are Parseval frame vectors for $\H$ with
 $p'_{i} = p'_{\xi_{i}}$, and $\xi_{N+1}$ has the property that
$\{\pi(g)\xi_{N+1}\,|\,g\in G\}$ is a Parseval frame for its
closed linear span. Since  $p'_{1} , ..., p'_{N}, p'_{N+1}$ are
mutually orthogonal, we get  that $\xi_1,\dots,\xi_{N+1}$ is a
strongly disjoint $N+1$-tuple of Parseval frame vectors. Moreover,
from $\sum_{i=1}^{N+1}p'_{i} = I$ we obtain that
$\xi_1\oplus\dots\oplus\xi_N\oplus\xi_{N+1}$ has norm one and
hence it is an ONB vector for $\pi^{N}\oplus \pi^{r}$, where
$\pi^{r} = \pi p'_{r}$ and $p'_{r}$ is the orthogonal projection
onto $\overline{span}\{\pi(g)\xi_{N+1}\,|\,g\in G\}$.

 (ii) If $\xi_{N+2}$ is strongly disjoint from all $\xi_i$,
$i=1,N$, then $\zeta:=(\xi_1\oplus+\dots\oplus
\xi_{N+1})\oplus\xi_{N+2}$ is a Parseval frame vector with
$\|\zeta\|>1$, and this is impossible.

(iii) This clearly is a more general form of $(i)$, and the proof
is exactly as in the proof of (i).
\end{proof}

The following two more lemmas are needed in the proof of Theorem
\ref{pri6}.

\begin{lemma}\cite{HaLa08} \label{pri4}
Let $G$ be a countable group and let $\pi:G\rightarrow\U(\H)$ be a
unitary representation of the group $G$ on the Hilbert space $\H$.
Assume there exists a Parseval frame vector $\xi_1$ for $\H$.
Suppose $\eta_1,\eta_2$ are two Parseval frame vectors for some
subspaces of $\H$. Let $p_{\eta_i}$ be the projection onto the
subspace $\cj{\pi(G)'\eta_i}$, $i=1,2$. Then

(i) the two Parseval frame vectors $\eta_1,\eta_2$ are strongly
disjoint if and only if projections $p_{\eta_1},p_{\eta_2}$ are
orthogonal.

(ii) the two Parseval frame vectors $\eta_1,\eta_2$ are unitarily
equivalent if and only if $p_{\eta_1}= p_{\eta_2}$.
\end{lemma}

\begin{lemma} \cite{GaHa01} \label{pri4-1} Let $G$ be a
countable group and let $\pi:G\rightarrow\U(\H)$ be a unitary
representation of the group $G$ on the Hilbert space $\H$. Assume
there exists a Parseval frame vector $\xi_1$ for $\H$. Suppose
$\eta$ is a Parseval frame vector for a subspace of $\H$. Then
there exists a vector $\zeta$ such that $\eta+\zeta$ is a Parseval
frame vector for $\H$, and $\eta$ and $\zeta$ are strongly
disjoint Parseval frame vectors.
\end{lemma}

\begin{proposition}\label{pri5}
Let $G$ be a countable ICC group and let $\pi:G\rightarrow\U(\H)$
be a unitary representation of the group $G$ on the Hilbert space
$\H$. Assume that there exists a Parseval frame vector $\xi_1$ for
$\H$. For $\xi\in\H$, let $p_\xi$ be the projection onto the
subspace $\cj{\pi(G)'\xi}$.
\begin{enumerate}
\item If $\eta$ is a Parseval frame vector for a subspace for $\H$
then $\tr_{\pi(G)''}(p_\eta)=\|\eta\|^2$. \item If $\eta$ is a
Parseval frame vector for $\H$ and $u$ is a unitary operator in
$\pi(G)''$, then the Parseval frame vector $u\eta$ is strongly
disjoint from $\eta$ iff $u\eta$ is orthogonal to the range of
$p_{\eta}$, and in this case $p_{\eta}\perp p_{u\eta}$. \item Let
$N\geq1, N\in\bz$, and suppose $\xi_1,\dots,\xi_{N+1}$ is a
strongly disjoint $N+1$-tuple of Parseval frame vectors for some
subspaces of $\H$, with $\sum_{i=1}^{N+1}\|\xi_i\|^2=1$. Then the
projections $p_{\xi_i}, i=1, \cdots , N$ are mutually orthogonal,
and $p_{\xi_1}+\dots+p_{\xi_{N+1}}=1$.
\end{enumerate}
\end{proposition}

\begin{proof}
By Proposition \ref{propi1}, we can assume $\H=p'l^2(G)$ and
$\pi=\lambda p'$ for some projection $p'\in L(G)'$ with
$\tr_{L(G)'}(p')=\|\xi_1\|^2$.

Let $p_{\eta}'$ be the projection onto the subspace
$\cj{\pi(G)''\eta}$ generated by the Parseval frame vector $\eta$.
Then $p_{\eta}\in\pi(G)'$ and by \cite[Remark 2.2.5]{JoSu97} we
have
$$\dim_{\pi(G)''}\H=\frac{\tr_{\pi(G)''}(p_{\eta})}{\tr_{\pi(G)'}(p_{\eta}')}.$$
But, by \cite[Proposition 2.2.6(vi)]{JoSu97},
$$\dim_{\pi(G)''}\H=\dim_{L(G)p'}(p'l^2(G))=\tr_{L(G)'}(p')\dim_{L(G)}l^2(G)=\|\xi_1\|^2.$$

On the other hand, using the uniqueness of the trace on factors,
we
have$$\tr_{\pi(G)'}(p_{\eta}')=\tr_{p'L(G)'p'}(p_{\eta}')=\frac{\tr_{L(G)'}(p_{\eta}')}{\tr_{L(G)'}(p')}=\frac{\|\eta\|^2}{\|\xi_1\|^2}.$$
We used the fact that $p_{\eta}'\in L(G)'$ is the projection onto
the subspace in $\H\subset l^2(G)$ spanned by the Parseval frame
$\{\pi(g)\eta=\lambda(g)\eta\,|\,g\in G\}$, and therefore, by
Proposition \ref{propi1}, $\tr_{L(G)'}(p_{\eta}')=\|\eta\|^2$. From these equalities, (i) follows.

To prove (ii), we use Theorem \ref{thi2} to see that $u\eta$ is a
Parseval frame vector. By Lemma \ref{pri4}, the Parseval frame
vectors are disjoint iff $p_{\eta}$ and $p_{u\eta}$ are
orthogonal. Thus one implication is trivial. If $u\eta$ is
orthogonal to the range of $p_{\eta}$, then for all
$x',y'\in\pi(G)'$ we have
$\ip{x'u\eta}{y'\eta}=\ip{u\eta}{x'^*y'\eta}=0$ so $p_{u\eta}$ is
perpendicular to $p_{\eta}$.

For (iii), with Lemma \ref{pri4}, we have that $p_{\xi_i}$ are
mutually orthogonal. From (i), we have that
$$\tr_{\pi(G)''}(\sum_{i=1}^{N+1}p_i)=\sum_{i=1}^{N+1}\|\xi_i\|^2=1,$$
therefore $\sum_{i=1}^{N+1}p_i=1$.

\end{proof}

Now we are ready to parametrize the Parseval frame vectors for $\H$.
As we mentioned before, such a result was given in \cite{HaLa00},
see Theorem \ref{thi2}. The drawback of their result is that they
are using operators in the von Neumann algebra $\pi(G)''$ itself,
not in its commutant. We give here a parametrization that uses
operators in the commutant. As Han and Larson proved in their
paper, we cannot expect to use unitary operators in $\pi(G)'$.
Instead, we will use the $N+1$ strongly disjoint Parseval frame
vectors given by Theorem \ref{thi3} and $N+1$ operators in the
commutant $\pi(G)'$ that satisfy the orthogonality relation
\eqref{eqi5_2} which in fact represents the norm-one property of a
first row in a unitary matrix of operators.

\begin{theorem}\label{pri6}
Let $G$ be a countable ICC group and let $\pi:G\rightarrow\U(\H)$
be a unitary representation of the group $G$ on the Hilbert space
$\H$. Suppose $\xi_i$, $i=1, \cdots, N+1$ is a strongly disjoint
$N+1$-tuple of Parseval frame vectors as in Theorem \ref{thi3}(i).
Let $\eta\in\H$ be another Parseval frame vector for a subspace of
$\H$, and let $q'$ be the projection onto this subspace. Then
there exist unique $u_i'\in \pi(G)'$, $i=1, \cdots, N+1$, with
$q'u_i'=u_i'$, $i=1, \cdots, N+1$, and $u_{N+1}'p_r'=u_{N+1}'$
(where $p_r'$ is the projection onto the span of the Parseval
frame generated by the vector $\xi_{N+1}$, as in Theorem
\ref{thi3}) such that
\begin{equation}\label{eqi5_1}
\eta=u_1'\xi_1+\dots+u_{N+1}'\xi_{N+1}.
\end{equation}
 Moreover
 \begin{equation}\label{eqi5_2}
 \sum_{i=1}^{N+1}u_i'u_i'^*=q'.
 \end{equation}
 Conversely, if the vector $\eta$ is defined by \eqref{eqi5_1} with $u_i'\in L(G)'$, $q'u_i'=u_i'$, $i=1,\cdots, N+1$, and $u_{N+1}'p_r'=u_{N+1}'$ satisfying \eqref{eqi5_2}, then $\eta$ is a Parseval frame vector for $\H$.
\end{theorem}

\begin{proof}
We prove the theorem first for the case when $q'=1$, so $\eta$ is
a Parseval frame vector for the entire space $\H$.

By Theorem \ref{thi3}(iii), there exist
$\eta_2,\dots,\eta_{N+1}\in \H$ that together with $\eta_1:=\eta$
form a strongly disjoint $N+1$-tuple of Parseval frame vectors as
in Theorem \ref{thi3}(i). Then $\eta_1\oplus\dots\oplus\eta_{N+1}$
is an ONB vector for $\pi^N\oplus \pi p_r'=\lambda$. Then there
exists a unitary $u'\in L(G)'$ such that
$u'(\xi_1\oplus\dots\oplus\xi_{N+1})=\eta_1\oplus\dots\oplus\eta_{N+1}$.

Let $p_i'$ be the projections onto the the $i$-th component. We
can identify $p_{N+1}'=p_r'$, and in our case $q'=p_1'$. Let
$u_i':=p_1'u'p_i'$. Then $u_i'\in\pi(G)'$,
$u_{N+1}'p_r'=u_{N+1}'$, and
$$\sum_{i=1}^{N+1}u_i'\xi_i=\sum_{i=1}^{N+1}p_1'u'p_i'(\xi_1\oplus\dots\oplus\xi_{N+1})=p_1'u'(\sum_{i=1}^{N+1}p_i')(\xi_1\oplus\dots\oplus\xi_{N+1})=\eta_1.$$
This proves \eqref{eqi5_1}.

To prove uniqueness, suppose $\sum_{i=1}^{N+1}v_i'\xi_i=0$ for
some operators $v_i'\in \pi(G)'$, with $v_{N+1}'p_r'=v_{N+1}'$.
Then, for all $g\in G$, $\pi(g)\sum_{i=1}^{N+1}v_i'\xi_i=0$. By
Lemma \ref{pri4}, the vectors $v_i'\xi_i$ are mutually orthogonal.
Since $\pi(g)$ is unitary, it follows that $\pi(g)v_i'\xi_i=0$ for
all $i=1, \cdots, N+1$. Therefore $v_i'\pi(g)\xi_i=0$ for all
$g\in G$, $i=1,N+1$. But $\pi(g)\xi_i$ span $\H$ for $i=1, \cdots,
N$, and span $p_r'\H$ for $i=N+1$. Therefore $v_i'=0$ for all
$i=1, \cdots, N+1$. This implies the uniqueness.

We check now \eqref{eqi5_2}. We have
$$q'=p_1'=p_1'u'u'^*p_1'=p_1'u'(\sum_{i=1}^{N+1}p_i)u'^*p_1'=\sum_{i=1}^{N+1}p_1'u'p_i'(p_1'u'p_i')^*=\sum_{i=1}^{N+1}u_i'u_i'^*.$$
This proves \eqref{eqi5_2}.

For the converse, we can use \cite[Proposition 2.21]{HaLa00}. We
include the details. Consider the frame transforms $\Theta_i$ for
the Parseval frame vectors $\xi_i$, $\Theta_i$ defined on $\H$ for
$i=1, \cdots, N$ and on $p_r'\H$ for $i=N+1$. We have, by
Proposition \ref{propi1}, $\Theta_i^*\Theta_i=1_\H$ for $i=1, ...
,N$, $\Theta_{N+1}^*\Theta_{N+1}=1_{p_r'\H}$. Note that
$\Theta_{u_{i}'\xi_{i}} = \Theta_{\xi_{i}}u_{i}'^{*}$ for $i =1,
... , N+1$ and $\Theta_{i}^{*}\Theta_{j} = 0$ for $i\neq j$. So we
have
$$
\Theta_{\eta}^{*}\Theta_{\eta}
=(\sum_{i=1}^{N+1}u_i'\Theta_i^*)(\sum_{i=1}^{N+1}\Theta_{i}u_i'^{*})
=\sum_{i,j=1}^{N+1}u_i'\Theta_i^*\Theta_ju_j'^*=\sum_{i=1}^{N+1}u_i'u_i'^*=I.
$$
Hence $\eta$ is a Parseval frame.

In the case when $q'\leq1$, we have $q'\in \pi(G)'$ and, by Lemma
\ref{pri4-1}, there exists a Parseval frame vector $\tilde\eta$
for $\H$ such that $\eta=q'\tilde\eta$. Using the proof above we
can find $\tilde u_i'\in\pi(G)'$ such that
$\tilde\eta=\sum_{i=1}^{N+1}\tilde u_i'\xi_i$ and all the other
properties. Then we can define $u_i':=q'\tilde u_i'$, $i=1,
\cdots, N+1$, and a simple computation shows that the required
properties are satisfied. The proof of uniqueness and the converse
in this case is analogous to the one provided for the case $q'=1$.
\end{proof}

We can use the above parametrization result to characterize
strongly disjoint  (resp.  unitary equivalent) Parseval frame
vectors in terms of the given parametrization. Recall that two
Parseval frame vectors $\eta$ and $\xi$  are unitary equivalent if
and only if their analysis operators have the same range spaces,
which in turn is equivalent to the condition that
$\Theta_{\eta}\Theta_{\eta}^{*} = \Theta_{\xi}\Theta_{\xi}^{*}$.

\begin{theorem}\label{pri7}
Let $G$ be a countable ICC group and let $\pi:G\rightarrow\U(\H)$
be a unitary representation of the group $G$ on the Hilbert space
$\H$. Suppose $\xi_i$, $i=1, \cdots, N+1$ is a strongly disjoint
$N+1$-tuple of Parseval frame vectors as in Theorem \ref{thi3}(i).
Let $\eta,\zeta\in\H,$ be two Parseval frame vectors for some
subspaces of $\H$. Suppose
$$\eta=u_1'\xi_1+\dots+u_{N+1}'\xi_{N+1},\quad \zeta=v_1'\xi_1+\dots+v_{N+1}'\xi_{N+1},$$
with $u_i',v_i'\in \pi(G)'$, $i=1,N+1$, $u_{N+1}'p_r'=u_{N+1}'$,
$v_{N+1}'p_r'=v_{N+1}'$. Then we have

(i)  $\eta$ and $\zeta$ are strongly disjoint if and only if
\begin{equation}\label{eqi5.1}
v_1'u_1'^*+\dots v_{N+1}'u_{N+1}'^*=0.
\end{equation}
(ii) $\eta$ and $\zeta$ are unitary equivalent if and only if
\begin{equation}\label{eqi5.2}
[u_1'^*, ... , u_{N+1}'^*]^{t} [u_1', ... , u_{N+1}']=[v_1'^*, ...
, v_{N+1}'^*]^{t} [v_1', ... , v_{N+1}'],
\end{equation}
where ``t'' represents the transpose of the row vector.
\end{theorem}

\begin{proof} (i)
Let $\psi,\psi'\in \H$. Then, by Proposition \ref{pr4}, $\eta$ and
$\xi$ are disjoint if and only if
$$0=\sum_{g\in G}\ip{\pi(g)\eta}{\psi}\cj{\ip{\pi(g)\zeta}{\psi'}}=\sum_{g\in G}\ip{\pi(g)(\sum_{i=1}^{N+1}u_i'\xi_i)}{\psi}\cj{\ip{\pi(g)(\sum_{j=1}^{N+1}v_j'\xi_j)}{\psi'}}=$$
$$\sum_{i,j=1}^{N+1}\sum_{g\in G}\ip{\pi(g)\xi_i}{u_i'^*\psi}\cj{\ip{\pi(g)\xi_j}{v_j'^*\psi'}}=(\mbox{since $\xi_i$ are mutually disjoint})=$$
$$\sum_{i=1}^{N+1}\sum_{g\in G}\ip{\pi(g)\xi_i}{u_i'^*\psi}\cj{\ip{\pi(g)\xi_i}{v_i'^*\psi'}}=$$
$$(\mbox{since $\xi_i$ is a Parseval frame vector}, u_{N+1}'^*\psi,v_{N+1}'^*\psi\in p_r'\H)=$$
$$\sum_{i=1}^{N+1}\ip{u_i'^*\psi}{v_i'^*\psi}=\ip{\sum_{i=1}^{N+1}v_i'u_i'^* \psi}{\psi'}.$$
Since $\psi,\psi'\in\H$ are arbitrary, the proof of (i) is
complete.

(ii)  Let $\Theta_{i}$ be the analysis operator for $\xi_{i}$.
Then $\Theta_{\eta} = \sum_{i=1}^{N+1}\Theta_{i}u_{i}'^* $ and
$\Theta_{\zeta} = \sum_{i=1}^{N+1}\Theta_{i}v_{i}'^* $. Then
$\eta$ and $\zeta$ are unitary equivalent if and only if $
\Theta_{\eta}\Theta_{\eta}^{*} = \Theta_{\zeta}\Theta_{\zeta}^{*},
$ i.e.
$$
\sum_{i=1}^{N+1} \Theta_{i}u_{i}'^*\Theta_{\eta}^* =
\sum_{i=1}^{N+1} \Theta_{i}v_{i}'^*\Theta_{\zeta}^*.
$$
Since $\Theta_{i}$ have orthogonal range spaces, we have that the
above equation holds if and only if
\begin{equation}\label{eqi5.3}
\Theta_{i}u_{i}'^*\Theta_{\eta}^* =
\Theta_{i}v_{i}'^*\Theta_{\zeta}^*
\end{equation} for all $i=1,
..., N+1$ . Applying $\Theta_{i}^{*}$ to both sides of
(\ref{eqi5.3}) and using that fact that $\Theta_{i}^{*}\Theta_{i}
= I$ for $i=1, ... , N$, $\Theta_{N+1}^{*}\Theta_{N+1} = p_{r}'$,
 $u_{N+1}'p_r'=u_{N+1}'$ and
$v_{N+1}'p_r'=v_{N+1}'$, we obtain that $ u_{i}'^*\Theta_{\eta}^*
= v_{i}'^*\Theta_{\zeta}^* $ for all $i=1, ..., N+1$, i.e.,
$$
\sum_{i=1}^{N}(u_{i}'^*u_{i}' - v_{i}'^* v_{i}')\Theta_{i}^{*} =
0.
$$
Apply the above left-side operator to $\Theta_{j}(\H)$ and use the
fact $\Theta_{i}^{*}\Theta_{j} = 0 $ when $i\neq j$, we get
$u_{i}'^*u_{i}' - v_{i}'^* v_{i}' = 0$ for all $i, j$. Hence we
have
$$
u_1'^*, ... , u_{N+1}'^*]^{t} [u_1', ... , u_{N+1}']=[v_1'^*, ...
, v_{N+1}'^*]^{t} [v_1', ... , v_{N+1}'].
$$
Conversely, if the above identity holds, then we clearly have $$
\sum_{i=1}^{N}(u_{i}'^*u_{i}' - v_{i}'^* v_{i}')\Theta_{i}^{*} = 0
$$
for all $i$ and so we have that (\ref{eqi5.3}) holds for all $i$.
Therefore $\eta$ and $\zeta$ are unnitary equivalent
\end{proof}


We conclude this section by pointing out that we also have similar
results as Theorem \ref{pri6} and Theorem \ref{pri7} for frame
representations of arbitrary countable groups. The following lemma
replaces Theorem \ref{thi3} for the general group case. Note that,
unlike the ICC group case, here we can not require that $\xi_{2},
... , \xi_{N}$ are Parseval frame vectors for $\H$. Recall that
the cyclic multiplicity for an subspace $\mathcal S$ of operators
on $\H$ is the smallest cardinality $k$ such that there exist
vectors $y_{i}$ ($i=1, ..., k)$ with the property
$\overline{span}\{Sy_{i}: S\in\mathcal S, i=1, ..., k\} = \H$.

\begin{lemma} \cite{HaLa08} \label{lemma5.7}
Let $\pi:G\rightarrow\U(\H)$ be a unitary representation of a
countable group $G$ on the Hilbert space $\H$ such that $\pi(G)'$
has cyclic multiplicity $N+1$ (here $N$ could be $\infty$). Assume
$\xi_{1}$ is a Parseval frame vector for $\H$. Then there exist
$\xi_{i}$ for $i=2, ... , N+1$ with the properties:

(i)\ $\{\pi(g)\xi_{i}: g\in G\}$ is a Parseval frame for $M_{i}: =
\overline{span}\{\pi(g)\xi_{i}: g\in G\}$;

(ii) \ $\xi_{i}$ ($i=1, ... , N+1$) are mutually strongly
disjoint;

(iii) there is no non-zero Bessel vector which is strongly
disjoint with all $\xi_{i}$.

\end{lemma}

One fact we used in the proof of Theorem \ref{pri6} is that two
ONB vectors for a unitary representation are linked by a unitary
operator in the commutant of the representation, i.e. all the ONB
vectors are unitarily equivalent. However, as we have already
mentioned before this is no longer true in general for
non-ONB vectors. The following characterizes the representations
that have this property, and it is needed in order to prove our
new parametrization result (Theorem \ref{thm-general}) for general
countable groups.

\begin{lemma}\cite{GaHa04}  \label{lemma5.8}
Let $\pi:G\rightarrow\U(\H)$ be a unitary representation of a
countable group $G$ on the Hilbert space $\H$ and $\xi$ be a
Parseval frame vector for $\H$. Then the following statements are
equivalent

(i) there is no non-zero Bessel vector which is strongly disjoint
with all $\xi$;

(ii) $P_{\xi} = \Theta_{\xi}\Theta_{\xi}^{*}\in L(G)'\cap L(G)''$;

(iii) a vector $\eta$ is a Parserval vector for $\H$ if and only
if there exists a unitary $u'\in \pi(G)'$ such that $\eta =
u'\xi$.
\end{lemma}

\begin{theorem}\label{thm-general}
Let $G$ be a countable group and let $\pi:G\rightarrow\U(\H)$ be a
unitary representation of the group $G$ on the Hilbert space $\H$.
Suppose $\xi_i$, $i=1,..., N+1$ are as in Lemma \ref{lemma5.7}.

(i) Let $\eta\in\H$. Then $\eta$ is a Parseval frame vector for
$\H$ if and only if there exist $u_i'\in \pi(G)'$, $i=1,..., N+1$,
with $u_i'p'_i=u_i'$, $i=1,... , N+1$ (where $p_i'$ is the
projection onto the closed linear span of $\{\pi(g)\xi_{i}:g\in
G\}$) such that $$ \eta=u_1'\xi_1+\dots+u_{N+1}'\xi_{N+1}$$ and
$\sum_{i=1}^{N+1}u_{i}'u_{i}'^* = I$.  Moreover these $u_{i}'$s
are unique.

(ii) Let $\eta,\zeta\in\H,$ be two Parseval frame vectors for $\H$
such that
$$\eta=u_1'\xi_1+\dots+u_{N+1}'\xi_{N+1},\quad \zeta=v_1'\xi_1+\dots+v_{N+1}'\xi_{N+1},$$
with $u_{i}$ and $v_{i}$ satisfying the requirement as in (i).
Then

(a) $\eta$ and $\zeta$ are strongly disjoint if and only if $$
v_1'u_1'^*+\dots v_{N+1}'u_{N+1}'^*=0,$$ and

(b) $\eta$ and $\zeta$ are unitary equivalent if and only if $$
[u_1'^*, ... , u_{N+1}'^*]^{t} [u_1', ... , u_{N+1}']=[v_1'^*, ...
, v_{N+1}'^*]^{t} [v_1', ... , v_{N+1}']. $$

\end{theorem}
\begin{proof} We only give a sketch proof for the necessary part
of $(i)$. The rest is similar to the ICC group case.

Let $\eta$ is a Parseval frame vector for $\H$. Then we have that
$p'_{\eta}\sim p'_{\xi_{1}}$ in $L(G)'$ by Lemma \ref{lemma2.1}.
Since $L(G)'$ is a finite von Neumann algebra,
we can find projections \cite{KaRi97} $q_{i}'$$(i=2, ..., N+1)$
such that $p_{\eta}'$, $q_{i}'$$(i=2, ..., N+1)$ are mutually
orthogonal and $q_{i}'\sim p_{\xi_{i}}'$ for $i =2, ... , N+1$.
From Lemma \ref{lemma2.1}, there exist Parseval frame vectors
$\eta_{i} $ for $M_{i} := \overline{span}\{\pi(g)\xi_{i}: g\in
G\}$ such that $p'_{\eta_{i}} = q'_{i}$ $(i =2, ... , N+1)$.
Define unitary representation $\sigma: G \rightarrow \U(\K)$
(where $K = \H\oplus M_{2}\oplus \cdots \oplus M_{N+1}$) by
$$
\sigma(g) = \pi(g) \oplus \pi(g)p_{\xi_{2}}\oplus \cdots \oplus
\pi(g)p_{\xi_{N+1}}
$$
Then both $\tilde{\xi}: = \xi_{1}\oplus\cdots\oplus \xi_{N+1}$ and
$\tilde{\eta}: = \eta_{1}\oplus\cdots\oplus \eta_{N+1}$ are
Parserval frame vectors for $\K$. By Lemma \ref{lemma5.7} (iii) we
have that there is no non-zero Bessel vector that is strongly
disjoint with all $\xi_{i}$. This implies that there is no
non-zero Bessel vector that is strongly disjoint with
$\tilde{\xi}$. Thus, from  Lemma \ref{lemma5.8},  we get that
there is a (unique) unitary operator $u'\in \sigma(G)'$ such that
$\tilde{\eta} = u'\tilde{\xi}$. Let $u_{i}' = u_{1,i}'P_{i}'$
(where $[u'_{1,1}, ... , u'_{1, N+1}]$ is the first row vector of
$u'$. Then it can be checked that $
\eta=u_1'\xi_1+\dots+u_{N+1}'\xi_{N+1},  $ and $ u_{i}'$ satisfy
all the requirements listed in $(i)$.
\end{proof}

We make a final remark that all the results  in this section
remain valid when unitary representation is replaced by {\it
projective unitary representations} of countable groups. In
particular, Theorem \ref{pri6} and Theorem \ref{pri7} remain true
for Gabor unitary representations. The interested reader can check
(cf. \cite{GaHa01, GaHa03, GaHa04, Gro, Han3, Han4, HaLa08,
Heil07}) for definitions and recent developments about projective
unitary representations and Gabor representations.

\section{Parseval frames and subgroups}\label{pars}

In this section we will only be interested in Parseval frame
vectors $\xi$ that have $\|\xi\|^2=\frac{1}{N}$. We will assume in
addition that there is a subgroup $H$ of $G$ of index $N$. The
next proposition shows that in this situation we can find
orthonormal bases for $\H$ obtained by the action of the subgroup
$H$.

\begin{proposition}\label{prp1}
Let $G$ be a countable ICC group and let $\pi:G\rightarrow\U(\H)$
be a unitary representation of the group $G$ on the Hilbert space
$\H$. Suppose $\xi\in\H$ is a Parseval frame vector for $\H$ with
$\|\xi\|^2=\frac{1}{N}$ for some $N\in\bn$. Assume in addition
that there exists an ICC subgroup $H$ of index $[G:H]=N$. Then
there exists a Parseval frame vector $\eta$ for $\H$ with the
property that $\sqrt{N}\{\pi(h)\eta\,|\,h\in H\}$ is an
orthonormal basis for $\H$.
\end{proposition}

\begin{proof}
By Theorem \ref{thi2} we can assume that $\H=p'l^2(G)$ for some
$p'\in L(G)'$ with $\tr_{L(G)'}(p')=\frac1N$ and $\pi=\lambda$
restricted to $p'l^2(G)$. We claim that $\dim_{L(H)}\H=1$.

Using \cite[Proposition 2.3.5, Example 2.3.3, Proposition
2.2.1]{JoSu97} we have
$$\dim_{L(H)}\H=\dim_{L(G)}\H\cdot [L(G):L(H)]=N\dim_{L(G)}p'l^2(G)=N\tr_{L(G)'}(p')\dim_{L(G)}l^2(G)=1.$$
Thus (see \cite[Chapter 2.2]{JoSu97}) the Hilbert space $\H$
considered as a module over $L(H)$, with the representation
$\pi(h)=\lambda(h)p'$, $h\in H$, is isomorphic to the module
$l^2(H)$, i.e., there exists an isometric isomorphism
$\Phi:\H\rightarrow l^2(H)$ such that $\Phi\pi(h)=\lambda(h)\Phi$
for all $h\in H$.

Define $\eta:=\frac{1}{\sqrt{N}}\Phi^{-1}(\delta_e)$. Then
$\sqrt{N}\{\pi(h)\eta\,|\,h\in H\}=\Phi^{-1}\{\delta_h\,|\,h\in
H\}$ so it is an orthonormal basis for $\H$.

We check that $\{\pi(g)\eta\,|\,g\in G\}$ is a Parseval frame for
$\H$. Let $\{a_1,\dots, a_N\}$ be a complete set of representatives
for the left cosets $\{gH\,|\,g\in G\}$. Let $v\in\H$. Then, since
$\sqrt{N}\{\lambda(h)\eta\,|\,h\in H\}$ is an orthonormal basis,
we have
$$\sum_{g\in G}|\ip{v}{\pi(g)\eta}|^2=\sum_{i=1}^N\sum_{h\in H}|\ip{v}{\pi(a_ih)\eta}|^2=\sum_{i=1}^N\sum_{h\in H}|\ip{\pi(a_i)^*v}{\pi(h)\eta}|^2=$$
$$\sum_{i=1}^N\frac{1}{N}\|\pi(a_i)^*v\|^2=1.$$
\end{proof}

\begin{remark}\label{rem3_1}
The condition $\|\xi\|^2=\frac1N$ is essential. In other words,
suppose $\eta$ is a Parseval frame vector for the representation
$\pi$ of $G$ on $\H$. Let $H$ be a subgroup of $G$ of index $N$,
and suppose the family $\{\pi(h)\eta\,|\,h\in H\}$ is an
orthogonal basis for the {\it whole} space $\H$. Then
$\|\xi\|^2=\frac1N$.

To see this, let $g_1\dots g_N\in G$ be representatives of the
left-cosets of $H$ in $G$. We have on one hand, using the
orthogonal basis, for all $x\in\H$:
$$\sum_{h\in H}|\ip{\pi(h)\eta}{x}|^2=\|\eta\|^2\|x\|^2.$$
On the other hand, using the Parseval frame
$$\|x\|^2=\sum_{i=1}^N\sum_{h\in H}|\ip{\pi(g_i)\pi(h)\eta}{x}|^2=\sum_{i=1}^N|\ip{\pi(h)\eta}{\pi(g_i)^*x}|^2=$$
$$\sum_{i=1}^N\|\eta\|^2\|\pi(g_i)^*x\|^2=N\|x\|^2\|\eta\|^2.\mbox{  Thus }\|\eta\|^2=\frac{1}{N}.$$
\end{remark}

We saw in Theorem \ref{thi3} that we can construct $N$ strongly
disjoint Parseval frame vectors for our Hilbert space $\H$. We
want to see if we can do this in such a way that, by undersampling
with the subgroup $H$, we have orthonormal bases (up to a
multiplicative constant), as in Proposition \ref{prp1}. We prove
that this is possible in the case when $H$ is normal and has an
element of infinite order.

\begin{theorem}\label{thp2}
Let $G$ be a countable ICC group and let $\pi:G\rightarrow\U(\H)$
be a unitary representation of the group $G$ on the Hilbert space
$\H$. Suppose there exists a Parseval frame vector $\xi\in\H$ with
$\|\xi\|^2=\frac1N$, $N\in\bz$. Assume in addition that $H$ is a
normal ICC subgroup of $G$ with index $[G:H]=N$, and $H$ contains
elements of infinite order. Then there exist a strongly disjoint
$N$-tuple $\eta_1,\dots,\eta_N$ of Parseval frame vectors for $\H$
such that for all $i=1,N$, the family
$\sqrt{N}\{\pi(h)\eta_i\,|\,h\in H\}$ is an orthonormal basis for
$\H$.
\end{theorem}

\begin{proof}
By Proposition \ref{propi1}, we can assume that $\pi$ is the
restriction of the left regular representation $\lambda$ on
$p'l^2(G)$, where $p'$ is a projection in $L(G)'$, with
$\tr_{L(G)'}(p')=\frac1N$.

We will define some unitary operators $u_i$ on $l^2(G)$ that will
help us build the frame vectors $\eta_i$ from just one such frame
vector $\eta$ given by Proposition \ref{prp1}.

Let $a_k$, $k=0,N-1$ be a complete set of representatives for the
cosets in $G/H$. Since $H$ is normal, $a_k^{-1}H$, $k=0,N-1$ forms
a partition of $G$. We can take $a_0=e$.

Define the functions $\varphi_j:G\rightarrow\bc$,
$\varphi_j(g)=e^{2\pi i \frac{kj}{N}}$ if $g\in a_k^{-1}H$. Note
that
\begin{equation}
    \sum_{k=0}^{N-1}\varphi_i(a_k^{-1}g)\cj{\varphi_j(a_k^{-1}g)}=0,\quad(g\in G, i\neq j).
    \label{eqp1}
\end{equation}
Indeed, if $g=a_rh$ for some $r\in\{0,\dots, N-1\}$ and $h\in H$,
then $a_k^{-1}g$, $k=0,N-1$ will lie in different sets of the
partition $\{a_l^{-1}H\}_{l=0,N-1}$, because $H$ is normal. Then
$$\sum_{k=0}^{N-1}\varphi_i(a_k^{-1}g)\cj{\varphi_j(a_k^{-1}g)}=\sum_{k=0}^{N-1} e^{2\pi i\frac{(i-j)k}{N}}=0.$$
Also $\varphi_i(gh^{-1})=\varphi_i(g)$, for all $g\in G,h\in H$,
and, since $H$ is normal $\varphi_i(hg)=\varphi_i(g)$, $i=0,N-1$,
$g\in G$, $h\in H$.

Define the operators $u_j$ on $l^2(G)$ by
$u_j\delta_g:=\varphi_j(g)\delta_g$, for all $g\in G$, $j=0,N-1$.
Since $u_j$ maps an ONB to an ONB, it is a unitary operator on
$l^2(G)$.

Then for all $g\in G$, using \eqref{eqp1},
$$\sum_{k=0}^{N-1}\lambda(a_k)u_iu_j^*\lambda(a_k)^*\delta_g=\sum_{k=0}^{N-1}\varphi_i(a_k^{-1}g)\cj{\varphi_j(a_k^{-1}g)}\delta_{g}=0.$$

Also
$$
    u_i\lambda(h)\delta_g=\varphi_i(hg)\delta_{gh}=\varphi_i(g)\delta_{gh}=\lambda(h)u_i\delta_g,
$$
so $u_i$ commutes with $\lambda(h)$ for all $i=0,N-1$, $h\in H$.

We want to compress the unitaries $u_i$ to a subspace $p_1'l^2(G)$
for some well chosen projection $p_1'$ in $L(G)'$. Take $h_0\in
H$, such that $h_0^n\neq e$ for all $n\in\bz\setminus\{0\}$. Then,
with $\rho$ the right-regular representation,
$$\ip{\rho(h_0)^n\delta_e}{\delta_e}=\delta_n=\int_{\bt}z^n\,d\mu(z),\quad(n\in\bz).$$
where $\mu$ is the Haar measure on $\bt$. Then let $p_1'$ be the
spectral projection $\chi_E(\rho(h_0))$, where $E$ is a subset of
$\bt$ of measure $\frac{1}{N}$. We have
$$\tr_{L(G)'}(p_1')=\ip{p_1'\delta_e}{\delta_e}=\ip{\chi_E(\rho(h_0))\delta_e}{\delta_e}=\int_{\bt}\chi_E\,d\mu=\frac{1}{N}.$$
Also, for $i=0,N-1$, $g\in G$,
$$u_i\rho(h_0)\delta_g=u_i\delta_{gh_0^{-1}}=\varphi_i(gh_0^{-1})\delta_{gh_0^{-1}}=\varphi_i(g)\delta_{gh_0^{-1}}=\rho(h_0)u_i\delta_g,$$
so $u_i\rho(h_0)=\rho(h_0)u_i$, and therefore $p_1'$ commutes with
all $u_i$, $i=0,N-1$. In addition, since $p_1'\in L(G)'$, it
commutes with $\lambda(g)$ for all $g\in G$.

Then we compute for $i\neq j$
$$\sum_{k=0}^{N-1}(\lambda(a_k)p_1')(p_1'u_ip_1')(p_1'u_jp_1')(\lambda(a_k)p_1')^*\delta_g=p_1'\sum_{k=0}^{N-1}\lambda(a_k)u_iu_j^*\lambda(a_k)^*\delta_g=0.$$

The operators $\tilde u_i:=p_1'u_ip_1'$ are unitary on
$p_1'l^2(G)$ because $p_1'$ commutes with $u_i$. Also,
$p_1'u_ip_1'$ commute with $\lambda(h)$ on $p_1'l^2(G)$.

We will couple these results with the following lemma to
finish the proof.
\begin{lemma}\label{lemp3}
Let $G$ be an ICC group with an ICC subgroup $H$ of index
$[G:H]=N$. Let $p_1'$ be a projection in $L(G)'$ with
$\tr_{L(G)'}(p_1')=\frac{1}{N}$ and let $\eta\in p_1'l^2(G)$ such
that $\sqrt{N}\{\lambda(h)\eta\,|\,h\in H\}$ is an orthonormal
basis for $p_1'l^2(G)$. Suppose $\tilde u_i$, $i=0,N-1$ are
unitary operators on $p_1'l^2(G)$ such that $\tilde u_i$ commutes
with $\lambda(h)$ for all $h\in H$, and  for some complete set of
representatives $a_0,\dots,a_{N-1}$ of the left-cosets in $G/H$,
$$\sum_{k=0}^{N-1}\lambda(a_k)\tilde u_i\tilde u_j^*\lambda(a_k)^*=0,\quad(i\neq j).$$
Then the vectors $\tilde u_i\eta_0$, $i=0, \cdots, N-1$ have the
following properties:
\begin{enumerate}
\item $\sqrt{N}\{\lambda(h)\tilde u_i\eta_0\,|\,h\in H\}$ is an
orthonormal basis for $p_1'l^2(G)$ for all $i=0, \cdots, N-1$.
\item $\tilde u_0\eta_0,\dots,\tilde u_{N-1}\eta_0$ is a strongly
disjoint $N$-tuple of Parseval frame vectors.
\end{enumerate}
\end{lemma}
\begin{proof}
Since $\tilde u_i$ commutes with $\lambda(h)$ for all $h\in H$,
property (i) follows immediately from the hypothesis. This implies
also that $\tilde u_i\eta$ is a Parseval frame vector for
$p_1'l^2(G)$ (see the proof of Proposition \ref{prp1}).

To check the strong disjointness, let $\Theta_i$ be the frame
transform of the vector $\tilde u_i\eta_0$. Let
$\Theta_0^H:p_1'l^2(G)\rightarrow l^2(H)$ be the frame transform
for the $\frac{1}{\sqrt{N}}$-orthonormal basis
$\{\lambda(h)\eta\,|\,h\in H\}$. Then
${\Theta_0^H}^*\Theta_0^H=\frac{1}{N}I$. We have for $v\in
p_1'l^2(G)$:

$$\Theta_i(v)=(\ip{v}{\lambda(g)\tilde u_i\eta})_{g\in G}=\left((\ip{v}{\lambda(a_k)\lambda(h)\tilde u_i\eta})_{h\in H}\right)_{k=0,N-1}=$$
$$\left((\ip{\tilde u_i^*\lambda(a_k)^*v}{\lambda(h)\eta_0})_{h\in H}\right)_{k=0,N-1}.\mbox{ Then for }v,v'\in p_1'l^2(G),$$
$$\ip{\Theta_i(v)}{\Theta_j(v')}=\sum_{k=0}^{N-1}\sum_{h\in H}\ip{\tilde u_i^*\lambda(a_k)^*v}{\lambda(h)\eta}\cj{\ip{\tilde u_j^*\lambda(a_k)^*v}{\lambda(h)\eta}}=$$
$$\sum_{k=0}^{N-1}\ip{\Theta_0^H(\tilde u_i^*\lambda(a_k)^*v)}{\Theta_0^H(\tilde u_j^*\lambda(a_k)^*v')}=\sum_{k=0}^{N-1}\ip{\lambda(a_k)\tilde u_j{\Theta_0^H}^*\Theta_0^H\tilde u_i^*\lambda(a_k)^*v}{v'}=0.$$
This proves that the frames are strongly disjoint.
\end{proof}

Returning to the proof of the theorem, we see that we can apply
Lemma \ref{lemp3}. Let $\eta$ be a Parseval frame vector in $p_1'l^2(G)$
such that $\sqrt{N}\{\lambda(h)\eta\,|\,h\in H\}$ is an ONB for
$p_1'l^2(G)$. It can be obtained from Proposition \ref{prp1}.
Then, using Lemma \ref{lemp3}, we get that
$\eta_i:=p_1'u_ip_1'\eta$ form a strongly disjoint $N$-tuple of
Parseval frames and $\sqrt{N}\{\lambda(h)\eta_i\,|\,h\in H\}$ are
ONBs for $p_1'l^2(G)$. Then we can move everything onto our space
$\H$, because $\tr_{L(G)'}(p_1')=\tr_{L(G)'}(p')$ so the
projections $p_1'$ and $p'$ are equivalent in $L(G)'$ and the
representations $\lambda$ on $p_1'l^2(G)$ and $\pi$ on $\H$ are
equivalent.
\end{proof}

\begin{remark} \label{remark-final} (i) There exist ICC groups with ICC subgroups of
any finite index. For example, let $F_N$ the free group on $N$ generators and $p\in\mathbb{N}$.
Also, let $\phi:F_N\rightarrow\mathbb{Z}_p$ a surjective group morphism. Then $H:=\mbox{Ker} \phi$ is
a (free, thus ICC) normal subgroup of $F_N$ of finite index.\\
\\
(ii) There exist ICC groups without finite index proper subgroups.
For example, let $F$ be the Thompson's group and $F'$ its
commutator. Both groups are ICC (see e.g. \cite{Jol98}). Moreover,
$F'$ is infinite and simple. By a classic group theoretical
argument an infinite simple group cannot have finite index proper
subgroups. Indeed, let $G$ be infinite simple and $H$ of finite
index $k$ in $G$.  Then there exists a group morphism from $G$ to
the (finite) group of permutations of the set of right cosets
$X:=\{ Hg_i\mbox{}|\mbox{}i=1,\cdots, k\}$. This is given by
$g\rightarrow\alpha_g$ where  $\alpha_g(Hg_i)=Hg_ig$. Because $G$
is infinite the kernel of the above morphism must be non-trivial.
Moreover, because $G$ is simple the kernel
must be all of $G$, i.e. $\alpha_{g_i}(H)=H=Hg_i$ for all $i=1,\cdots, k$. Hence $H$ is not proper. \\
\\
(iii) If $G$ is ICC and $H$ is a finite index subgroup of $G$ then
$H$ is ICC. Indeed, let $G=\cup_{j=1}^{k} c_j H$, where $c_j$ are
distinct left cosets representatives. If for some $h\in H$ the
conjugacy class $\{ghg^{-1}\mbox{}|\mbox{}g\in H \}$ is finite
then the set $\{c_j g h (c_j g)^{-1}\mbox{}|\mbox{}g\in H,
j=1,\cdots, k\}$ is finite. Notice $\{c h
c^{-1}\mbox{}|\mbox{}c\in G\}\subset \{c_j g h (c_j
g)^{-1}\mbox{}|\mbox{}g\in H, j=1,\cdots ,k\}$.
However, the conjugacy class of $h$ in $G$ is infinite as $G$ is ICC. \\
\\
(iv) There exist ICC groups with all elements of finite order,
e.g. the Burnside groups of large enough exponents (see
\cite{Ol80}).

\end{remark}

\begin{proposition}\label{prp1.5}
Let $G$ be a countable ICC group and let $\pi:G\rightarrow\U(\H)$
be a unitary representation of the group $G$ on the Hilbert space
$\H$. Suppose $\xi\in\H$ is a Parseval frame vector for $\H$ with
$\|\xi\|^2=\frac{M}{N}$ for some $M,N\in\bn$. Assume in addition
that there exists a normal ICC subgroup $H$ of index $[G:H]=N$
such that $H$ contains elements of infinite order. Then there
exists $K:=\left\lfloor \frac NM\right\rfloor$ strongly disjoint
Parseval frame vectors $\eta_i$ for $\H$, $i=1, \cdots, K$ such
that $\{\pi(h)\eta_i\,|\,h\in H\}$ is an orthogonal family (in a
subspace of $\H$) for all $i=1, \cdots, K$.
\end{proposition}

\begin{proof}
Consider a projection $p'$ in $L(G)'$ of trace
$\tr_{L(G)'}(p')=\frac1N$. Using Theorem \ref{thp2} we can find
$\xi_1,\dots,\xi_N$ strongly disjoint Parseval frame vectors for
$\H_{1/N}:=p' l^2(G)$ with the representation
$\pi_{1/N}:=p'\lambda$, such that for all $i=1, \cdots, N$,
$\{\pi(h)\xi_i\,|\,h\in H\}$ is an orthogonal basis for
$\H_{1/N}$.

Then consider the representation $\pi_{1/N}^M$ on $\H_{1/N}^M$
with the strongly disjoint Parseval frame vectors
$\eta_i:=\xi_{(M-1)i+1}\oplus\dots\oplus\xi_{Mi}$, $i=1, \cdots,
K$.

Using Proposition \ref{propi1} we have that $\pi_{1/N}^M$ is
equivalent to a subrepresentation of the left-regular
representation, corresponding to a projection of trace
$\|\xi_1\oplus\dots\oplus\xi_M\|^2=\frac MN$. But the same is true
for the representation $\pi$ on $\H$. Therefore the two
representations are equivalent and the vectors $\eta_i$ can be
mapped into $\H$ to obtain the conclusion.
\end{proof}

\begin{remark}\label{rem3_2}
Suppose the hypotheses of Theorem \ref{thp2} are satisfied, with
$N\geq2$. Then we can construct uncountably many inequivalent
Parseval frame vectors $\eta$ for $\H$ with the property that
$\{\pi{h}(\eta)\,|\, h\in H\}$ is a Riesz basis for $\H$.

To see this, use Theorem \ref{thp2} to obtain strongly disjoint
Parseval frame vectors $\eta_1,\dots,\eta_N$ for $\H$, such that
$\{\pi(h)\eta_i\,|\,h\in H\}$ is an orthogonal basis for $\H$, for
all $i=1,N$.

Then take $\alpha,\beta\in\bc$ with $|\alpha|^2+|\beta|^2=1$, and
$|\alpha|\neq|\beta|$. Using Theorem \ref{pri6}, we obtain that
$\eta_{\alpha,\beta}:=\alpha\eta_1+\beta\eta_2$ is a Parseval
frame vector for $\H$.

Since $\eta_1$ and $\eta_2$ generate orthogonal bases under the
action of $H$, there is a unitary $u\in \pi(H)'$ such that
$u\eta_1=\eta_2$. Then $\eta_{\alpha,\beta}=(\alpha+\beta
u)\eta_1$. Since $u$ is unitary and $|\alpha|\neq |\beta|$ it
follows that $\alpha+\beta u$ is invertible. Therefore
$\{\pi(h)\eta_{\alpha,\beta}\,|\, h\in
H\}=\{(\alpha+u\beta)\pi(h)\eta_1\,|\,h\in H\}$ is a Riesz basis
for $\H$.

It remains to see when two such vectors $\eta_{\alpha,\beta},
\eta_{\alpha',\beta'}$ are equivalent. Using Theorem \ref{pri7} we
see that this happens only if $|\alpha|=|\alpha'|$,
$|\beta|=|\beta'|$ and $\cj\alpha\beta=\cj\alpha'\beta'$, i.e.,
$(\alpha,\beta)=c(\alpha',\beta')$ for some $c\in\bc$ with
$|c|=1$. Since we can find uncountably many pairs $(\alpha,\beta)$
such that no two such pairs satisfy this condition, it follows
that we can construct uncountably many inequivalent Parseval frame
vectors $\eta_{\alpha,\beta}$ that satisfy the given conditions.
\end{remark}

Finally we discuss how our results fit in the recent effort on the
Feichtinger's frame decomposition conjecture. It was recently
discovered (in particular, by Pete Casazza and his collaborators)
 that the famous intractible 1959
Kadison-Singer Problem in C*-algebras is equivalent to fundamental
open problems in a dozen different areas of research in
mathematics and engineering (cf. \cite{Cas06, Cas2}).
Particularly, the KS-problem is equivalent to the Feichtinger's
problem which asks whether every bounded frame (i.e. the norms of
the  vectors in the frame sequence are bounded from below) can be
written as a finite union of Riesz sequences. Since this question
is intractible in general, much of the effort has been focused on
special classes of frames. One natural and interesting class to
consider is the class of frames obtained by group representations
[See open problems posted at the 2006 ``The Kadison-Singer
Problem" workshop]. Unfortunately, except for a very few cases
(e.g., Gabor frames associated with rational lattices \cite{Cas05}) very little
is known so far even for this special class. Particularly,  it is
unknown whether for every (frame) unitary representation we can
always find one frame vector which is ``Riesz sequence"
decomposable. Therefore the results obtained in this section
certainly addressed some aspects of the research effort in this
direction. In particular, we have the following as a consequence
of our main result.

\begin{proposition}
Let $G$ be a countable ICC group and  assume that there exists an
ICC subgroup $H$ of index $[G:H]=N$. Then

(i) If $p'\in L(G)'$ is any projection such that $tr_{L(G)'}(p')
\ge \frac1N $,  then there exist a Parseval frame vector $\eta$
for the subrepresentation $\pi:=L|_{p'}$ such that $\{\pi(g)\eta:
g\in G\}$ is a finite union of Riesz sequences.

(ii) For any ONB vector for the left regular representation $L$
and any $\alpha$ such that $1>\alpha \ge \frac1N $, there exists a
projection $p'\in L(G)'$ such that $tr_{L(G)'}(p') = \alpha$ and
$\{p'L(g)\psi: g\in G\}$ is a finite union of Riesz sequences.

\end{proposition}

\begin{proof} (i) Since $L(G)'$ is a factor von Neumann algebra,
there exists a subprojection $q'$ of $p'$ such that
$tr_{L(G)'}(q') = \frac1N $. By Proposition \ref{prp1}, there
exists a Parseval frame vector, say $\eta_{1}$, for the
representation $\pi|_{q'}$ such that $\{\sqrt{N}\pi(h)\eta_{1}:
h\in H\}$ is orthonormal. By Lemma \ref{pri4-1}, we can ``dilate"
$\eta_{1}$ to a Parseval frame vector $\eta$ for $\pi$. Let
$\eta_{2} = (p'-q')\eta$. Then for any sequence $\{c_{h}\}_{h\in
H}$ (finitely many of them are non-zero) we have
$$
||\sum_{h\in H}c_{h}\pi(h)\eta||^{2} = ||\sum_{h\in
H}c_{h}\pi(h)\eta_{1}||^{2} + ||\sum_{h\in
H}c_{h}\pi(h)\eta_{2}||^{2} \geq ||\sum_{h\in
H}c_{h}\pi(h)\eta_{1}||^{2}
$$
since $\pi(G)\eta_{1}$ and $\pi(G)\eta_{2}$ are orthogonal. Thus
$\{\pi(h)\eta: h\in H\}$ is a Riesz sequence as $\{\pi(h)\eta_{1}:
h\in H\}$ is Riesz.

(ii) Let $r'\in L(G)'$  be a projection such that $tr_{L(G)'}(r')
= \alpha$. Then by part (i) there exists a Parseval frame vector
$\xi$ such that $\{\sigma(h)\xi: h\in H\}$ is a Riesz sequence,
where $\sigma = L|_{r'}$. We will show that there exists a
projection $p'\in L(G)'$ such that $\{p'L(g)\psi: g\in G\}$ and
$\{\sigma(g)\xi: g\in G\}$ are unitarily equivalent, and this will
imply that $\{p'L(h)\psi: h\in H\}$ is a Riesz sequence.

In fact, again by Lemma \ref{pri4-1}, there exists ONB vector
$\tilde{\psi}$ for the left regular representation $L$ such that
$r'\tilde{\psi} = \xi$. Since both $\psi$ and $\tilde{\psi}$ are
ONB vectors  for $L$ we have that there exists a unitary operator
$u'\in L(G)'$ such that $\psi = u' \tilde{\psi}$. Let $p' =
u'r'u'^{*}$. Then $p'\in L(G)'$ is a projection such that
$tr_{L(G)'}(p') =tr_{L(G)'}(r')=\alpha$. Moreover,
$$
p'L(g)\psi = L(g)p'\psi = L(g)(u'r'u'^{*})u'\tilde{\psi} =
u'L(g)r'\tilde{\psi} = u'L(g)\xi = u'\sigma(g)\xi
$$
for all $g\in H$. Hence $\{p'L(g)\psi: g\in G\}$ and
$\{\sigma(g)\xi: g\in G\}$ are unitarily equivalent, as claimed
and so we completed the proof.
\end{proof}

There are some interesting special cases. For example, as we
mentioned in Remark \ref{remark-final} if $G$ is a free group with
more than one generator, then we can find $N_{k}\rightarrow
\infty$ such that there exist ICC subgroups $H_{k}$ having the
property $[G:H_{k}] = N_{k}$.  Thus we have the following
corollary which for the free group case answered affirmatively one
of two open problems posted by Deguang Han at the 2006 ``The
Kadison-Singer Problem" workshop.

\begin{corollary}
Let $G$ be a free group with more than one generator. Then

(i) For any non-zero projection $p'\in L(G)'$,  there exists a
Parseval frame vector $\eta$ for the subrepresentation
$\pi:=L|_{p'}$ such that $\{\pi(g)\eta: g\in G\}$ is a finite
union of Riesz sequences.

(ii) For any ONB vector for the left regular representation $L$,
and any $\alpha > 0$, there exists a projection $p'\in L(G)'$ such
that $tr_{L(G)'}(p') = \alpha$ and $\{p'L(g)\psi: g\in G\}$ is a
finite union of Riesz sequences. This sequence will be orthogonal
when $\alpha = \frac1N $ for some $N\in\N$.

\end{corollary}

\bibliographystyle{alpha}
\bibliography{grfr}

\end{document}